\documentclass[12pt,reqno,letterpaper]{amsart}
\usepackage[T1]{fontenc}
\usepackage[utf8]{inputenc}
\usepackage{lmodern}
\usepackage[margin=2.5cm]{geometry}
\usepackage{amsmath,amssymb,amsthm,mathtools,booktabs,array}
\usepackage{microtype}
\usepackage{tikz}
\tikzset{finite graph/.style={every node/.style={draw=none,fill=white,
inner sep=1.5pt,font=\small}}}
\usepackage{listings}
\usepackage[colorlinks=true,linkcolor=blue,citecolor=blue,urlcolor=blue,unicode]{hyperref}
\usepackage[switch,mathlines]{lineno}
\makeatletter
\@ifundefined{subjclassname@2020}{%
  \@namedef{subjclassname@2020}{\textup{2020} Mathematics Subject Classification}%
}{}
\makeatother

\numberwithin{equation}{section}
\newtheorem{theorem}{Theorem}[section]

\newtheorem{lemma}[theorem]{Lemma}
\newtheorem{corollary}[theorem]{Corollary}
\theoremstyle{remark}

\newcommand{\R}{\mathbb R}

\newif\ifshowlineno
\showlinenofalse
\allowdisplaybreaks[2]
\title[Finiteness of five-body Dziobek configurations]{Finiteness of Dziobek central configurations in the Newtonian five-body problem}
\author{Kefeng Liu}
\address[Kefeng Liu]{Mathematical Sciences Research Center, Chongqing University of Technology, Chongqing, 400054, China}
\email{kefengliu@cqut.edu.cn}
\author{Yangshanshan Liu}
\address[Yangshanshan Liu]{Mathematical Sciences Research Center, Chongqing University of Technology, Chongqing, 400054, China}
\email{lyss6133@126.com}
\date{\today}
\hypersetup{pdfauthor={Kefeng Liu and Yangshanshan Liu}}
\begin{document}

\begin{abstract}
We prove that, for any five fixed strictly positive masses, the Newtonian
central configurations of exact affine dimension three form finitely many
similarity classes.
The proof uses the Dziobek equations, two constant identities, and local
analysis of zeros and poles on a compact complex curve.
\end{abstract}
\maketitle

\nolinenumbers
\tableofcontents
\ifshowlineno\linenumbers\fi

\section{Introduction}
\label{sec:introduction}

Consider $n$ bodies of fixed masses $m_i>0$ at distinct positions
$x_i\in\R^d$. Write $x=(x_1,\ldots,x_n)$,
$r_{ij}=|x_i-x_j|$, $m_0=\sum_{i=1}^n m_i$, and
$c=m_0^{-1}\sum_{i=1}^n m_ix_i$ for the center of mass.
With the gravitational constant set to one, the Newtonian equations are
\[
 m_i\ddot x_i=\sum_{j\ne i}m_i m_j\frac{x_j-x_i}{r_{ij}^3}
 =\frac{\partial\mathcal U}{\partial x_i},
 \qquad 1\le i\le n.
\]
The positive Newtonian potential function and the moment of inertia
about the center of mass are, respectively,
\[
 \begin{aligned}
 \mathcal U(x)&=\sum_{i<j}\frac{m_i m_j}{r_{ij}},\\
 \mathcal I(x)&=\sum_{i=1}^n m_i|x_i-c|^2
 =\frac{1}{m_0}\sum_{i<j}m_i m_jr_{ij}^2.
 \end{aligned}
\]
A \emph{central configuration} is a collision-free arrangement for which
every acceleration points toward the center of mass with the same
proportionality factor. Thus, for some $\lambda>0$,
\begin{equation}\label{eq:central-configuration}
 \sum_{j\ne i}m_j\frac{x_i-x_j}{r_{ij}^3}
 =\lambda(x_i-c),\qquad 1\le i\le n.
\end{equation}
Since $\nabla_{x_i}\mathcal U=-\sum_{j\ne i}m_i m_j(x_i-x_j)/r_{ij}^3$
and $\nabla_{x_i}\mathcal I=2m_i(x_i-c)$, this is equivalent to
\[
 \nabla_x\mathcal U+\frac{\lambda}{2}\nabla_x\mathcal I=0.
\]
Here $\nabla_x$ differentiates with respect to all positions.
The scaling laws $\mathcal U(tx)=t^{-1}\mathcal U(x)$ and
$\mathcal I(tx)=t^2\mathcal I(x)$, differentiated at $t=1$, give
\begin{equation}\label{eq:cc-energy}
 \mathcal U=\lambda\mathcal I,\qquad
 \frac{\lambda}{m_0}=1\ \Longrightarrow\ \mathcal U=m_0\mathcal I.
\end{equation}
These configurations generate homothetic motions, in which the shape is
preserved while its size changes.
The finiteness problem asks whether fixed positive masses admit only
finitely many central configurations up to similarity. Throughout this
paper, configurations are considered up to translations, orthogonal
transformations (including reflections), and positive dilations, with
the body indices fixed.

The planar version of this question is Smale's sixth problem for the
twenty-first century; see \cite[Sec.~2.14]{Moe15}.
Hampton and Moeckel proved finiteness for every positive mass vector in
the planar four-body problem
\cite[Thm.~1]{HM06}.
For the planar five-body problem, Albouy and Kaloshin proved finiteness
outside a closed algebraic subset of codimension two in the positive
mass space \cite[Thm.~2]{AK12}.
In \cite[Thm.~1.2]{LL26}, we prove finiteness
for every positive mass vector in the planar five-body problem using
analysis of singular sequences, potential bounds, and local
analytic reductions.
We study spatial five-body configurations of exact
affine dimension three.

A configuration of $n$ bodies is called a \emph{Dziobek configuration}
when its affine span, the smallest affine subspace containing its
positions, has dimension $n-2$
\cite[Def.~2.12.1]{Moe15}.
The name honors Otto Dziobek, who studied the planar four-body case
\cite{Dzi00}.
Thus spatial nonplanar five-body configurations are precisely the
five-body Dziobek configurations. Moeckel proved that Dziobek central
configurations have finitely many similarity classes for generic
masses, meaning masses outside a proper algebraic subset of the space
of mass ratios \cite[Thm.~3]{Moe01}; see also
\cite[Prop.~2.14.1]{Moe15}.
For the spatial five-body problem, Hampton and Jensen gave explicit
polynomial conditions on the masses outside which finiteness holds
\cite[Thm.~1; Tab.~1]{HJ11}.
A recent preprint of Dias also studies generic finiteness for Dziobek
central configurations with homogeneous potentials
\cite[Thm.~5.11]{Dia26}.
These results leave finiteness unresolved for exceptional mass vectors.

\begin{theorem}[Main theorem]
\label{js:theorem}
For every five fixed strictly positive masses, the Newtonian central
configurations of exact affine dimension three form finitely many
similarity classes.
\end{theorem}

The masses remain fixed throughout the proof, including when they satisfy
the equalities arising in the local analysis.

We follow Moeckel's notation
\cite[Sec.~2.12, Eqs.~(2.53)--(2.54)]{Moe15}.
For a five-body configuration $x_1,\ldots,x_5\in\R^3$ of affine
dimension three, form the augmented configuration matrix
\[
 \hat X=\begin{pmatrix}
 1&\cdots&1\\
 x_1&\cdots&x_5
 \end{pmatrix}\in\mathbb R^{4\times5}.
\]
Its rank is four, so its kernel has dimension one. Choose a nonzero
vector $\Delta=(\Delta_1,\ldots,\Delta_5)^T\in\ker\hat X$.
Then $\ker\hat X=\mathbb R\Delta$, and $\hat X\Delta=0$ is equivalent to
\[
 \sum_i\Delta_i=0,\qquad \sum_i\Delta_ix_i=0.
\]
Thus $\Delta$ records the affine dependence among the positions and is
determined up to a common nonzero factor. A Gale transform is the list
of rows of a matrix whose columns form a basis of $\ker\hat X$; here
it consists of the scalars $\Delta_i$.

For a central configuration, normalize $\lambda/m_0=1$. The Dziobek
equations then have a real coefficient $\kappa\ne0$
\cite[Prop.~2.12.2, Eq.~(2.60)]{Moe15}.
Subsection~\ref{js:setup} justifies this normalization, proves $\kappa<0$,
and introduces the normalized Gale coordinates
$\delta_i=\sqrt{-\kappa}\Delta_i/m_i$.
The distances then satisfy $r_{ij}^{-3}=1-\delta_i\delta_j$.
Equations \eqref{js:angular} describe constrained critical points with cube roots
$u_{ij}=(1-\delta_i\delta_j)^{1/3}>0$ and a common row value $\mu$.
Every physical Dziobek central configuration, with real positions and
positive Euclidean distances, supplies a solution of this system.
Finiteness of this solution set therefore implies the theorem.

The proof assumes that this solution set is infinite and constructs an
algebraic curve, an irreducible complex algebraic set of dimension one.
The curve is represented on a compact Riemann surface, a compact complex
manifold of dimension one
\cite[Chap.~I, Thm.~2.1; Chap.~II, Sec.~3,
Def.~3.1]{Gri89}.
The constant potential identity and positivity on the original real arc
force the cube-root germs to belong to the function field of the projected
curve (Lemma~\ref{js:constants}). They therefore define
single-valued meromorphic functions on its compact model, meaning complex
analytic functions with only poles as possible singularities
\cite[Chap.~I, Sec.~3, Def.~3.1]{Gri89}.
The construction is given in Subsection~\ref{js:setup}. Two identities then hold:
\[
 \sum_{i<j}m_i m_j u_{ij}
 =\sum_{i<j}m_i m_j u_{ij}^{-2}
 =U_0>\sum_{i<j}m_i m_j.
\]
The strict inequality uses positivity and the nonzero Gale vector.
For physical configurations, the two sums are $\mathcal U$ and
$m_0\mathcal I$, respectively, as in \eqref{eq:cc-energy}.

These identities restrict the zeros and poles of the coordinate
functions. At points where all $\delta_i$ are finite, the graph of
vanishing cube roots identifies a row with nonzero denominators,
so $\mu$ is holomorphic (complex analytic).
At any other point, exactly three or four coordinates have the largest
pole order. Their leading equations show that $\mu$ vanishes there and
cannot be identically zero. Thus $\mu$ is globally holomorphic on the
compact connected Riemann surface and must be constant
\cite[Chap.~I, Sec.~3, Thm.~3.3]{Gri89}.
A nonconstant $\delta_i$ has a pole, so this constant must be zero,
contradicting the local equations.

Subsection~\ref{js:setup} gives the reduction and the constant identities.
Subsection~\ref{js:finite-points} treats finite Gale points;
Subsection~\ref{js:poles} treats poles. The final contradiction appears
in Subsection~\ref{js:compact-contradiction}.

\section{Proof of the main theorem}
\label{js:section}

\subsection{Algebraic reduction and constant identities}
\label{js:setup}

From now on take $n=5$ and $d=3$, fix $m_1,\ldots,m_5>0$, and put
\[
 W=\sum_{i<j}m_i m_j.
\]
Consider a central configuration satisfying
\eqref{eq:central-configuration} with exact affine dimension three.
We use the notation introduced above and normalize $\lambda/m_0=1$
by a positive dilation, so that \eqref{eq:cc-energy} applies.

Use the affine dependence vector $\Delta$ defined in the Introduction.
By \cite[Prop.~2.12.2, Eq.~(2.60)]{Moe15}, there is a real
$\kappa\ne0$ such that
\[
 m_i m_j(r_{ij}^{-3}-1)=\kappa\Delta_i\Delta_j\qquad(i\ne j).
\]
Put $z_i=\Delta_i/m_i$, as in \cite[Sec.~2.14]{Moe15}. Then
\[
 r_{ij}^{-3}=1+\kappa z_iz_j,
 \qquad
 \sum_{i<j}\Delta_i\Delta_jr_{ij}^2
 =-\left|\sum_i\Delta_ix_i\right|^2=0.
\]
If $\kappa>0$, monotonicity of $(1+t)^{-2/3}$ on $t>-1$ gives
\[
 z_iz_j(1+\kappa z_iz_j)^{-2/3}\le z_iz_j,
\]
strictly when $z_iz_j\ne0$. A nonzero affine dependence has at least
two nonzero entries, so summation yields
\[
 0=\sum_{i<j}\Delta_i\Delta_jr_{ij}^2
 <\sum_{i<j}m_i m_jz_iz_j
 =-\frac12\sum_i m_i^2z_i^2<0.
\]
This contradiction gives $\kappa<0$.
Define
\[
 \delta_i=\sqrt{-\kappa}\,\frac{\Delta_i}{m_i}
 =\sqrt{-\kappa}\,z_i.
\]
We call the variables $\delta_i$ the \emph{normalized
Gale coordinates}; the entries $m_i\delta_i=\sqrt{-\kappa}\Delta_i$
form a Gale transform of the physical configuration.
The vector $\delta=(\delta_1,\cdots,\delta_5)$ is nonzero and satisfies
\[
 \sum_i m_i \delta_i=0,\qquad
 \sum_i m_i \delta_ix_i=0,\qquad
 r_{ij}^{-3}=1-\delta_i \delta_j>0.
\]
Replacing $\Delta$ by $a\Delta$, with $a\in\mathbb R\setminus\{0\}$,
replaces $\kappa$ by $\kappa/a^2$ and $\delta$ by
$\operatorname{sgn}(a)\delta$. Thus $\delta$ is unique up to its overall sign.

Put $u_{ij}=r_{ij}^{-1}=(1-\delta_i \delta_j)^{1/3}>0$.
Pair indices are unordered: $u_{ji}=u_{ij}$; the same convention applies
to all pair coefficients introduced below. The affine identities give
\[
 \sum_jr_{ij}^2m_j\delta_j=\sum_j|x_j|^2m_j\delta_j,
\]
independently of $i$. Hence every physical configuration supplies a
nonzero real solution of the critical-point equations
\begin{equation}\label{js:angular}
 u_{ij}^3=1-\delta_i \delta_j,\qquad
 \sum_i m_i \delta_i=0,\qquad
 \sum_{j\ne i}m_j\delta_j u_{ij}^{-2}=\mu\quad(1\le i\le5).
\end{equation}
With $u_{ij}>0$, these are the constrained critical
equations of
\[
 \Phi(\delta)=\sum_{i<j}m_i m_j(1-\delta_i \delta_j)^{1/3}
\]
on $\sum_i m_i \delta_i=0$. Its partial derivatives are
\[
 \frac{\partial\Phi}{\partial\delta_i}
 =-\frac{m_i}{3}\sum_{j\ne i}m_j\delta_j
 (1-\delta_i\delta_j)^{-2/3}.
\]
Thus constrained criticality is equivalent to $\nabla\Phi$ being a
multiple of the normal vector $(m_1,\ldots,m_5)$, or equivalently to
the five row sums in \eqref{js:angular} having a common value $\mu$.
The multiplier may depend on the critical point.
Equations \eqref{js:angular} give necessary conditions for a physical configuration.
Define the solution set of these constrained critical equations by
\[
 \mathcal S=\bigl\{(\delta,u,\mu)\in
 (\R^5\setminus\{0\})\times\R_{>0}^{10}\times\R:
 \text{\eqref{js:angular} holds}\bigr\}.
\]
Its projection onto the $\delta$ coordinates is precisely the set of
nonzero critical points of $\Phi$ subject to $\sum_i m_i\delta_i=0$,
in the domain $1-\delta_i\delta_j>0$ for all $i<j$.
Multiplying the $i$th critical equation by
$\prod_{j\ne i}u_{ij}^{2}>0$ clears its denominators, and
$\delta\ne0$ is the polynomial inequality $\sum_i\delta_i^2>0$.
Thus $\mathcal S$ is semialgebraic, meaning that it is described using
finitely many polynomial equalities and inequalities
\cite[Def.~2.1.4]{BCR98}.
Let $\mathcal D_{\mathrm{phys}}$ denote the distance vectors
$(r_{ij})_{1\le i<j\le5}$ of physical DCCs normalized by
$\lambda/m_0=1$. The reduction above gives
\[
 \mathcal D_{\mathrm{phys}}\subseteq
 \bigl\{(u_{ij}^{-1})_{i<j}:(\delta,u,\mu)\in\mathcal S\bigr\}.
\]
The distances determine the Gram matrix of the vectors $x_i-x_5$,
whose entries are their pairwise inner products:
\[
 \langle x_i-x_5,x_j-x_5\rangle
 =\frac{r_{i5}^2+r_{j5}^2-r_{ij}^2}{2},\qquad 1\le i,j\le4,
\]
where $r_{ii}=0$. Two realizations with the same Gram matrix give the
same squared norm to every corresponding linear combination.
Matching the relative position vectors therefore defines an isometry
between their spans, which extends to an orthogonal transformation.
Restoring the translation shows that the distances determine the
normalized configuration up to isometry. Hence finiteness of
$\mathcal S$ implies finiteness of physical similarity classes.

By \cite[Def.~9.1.7 and Prop.~9.1.8]{BCR98},
$\mathcal S$ admits a Nash stratification, that is, a finite partition into
semialgebraically connected Nash submanifolds, which are real-analytic
semialgebraic manifolds
\cite[Def.~2.9.9]{BCR98}. A stratum of dimension zero is a
single point, so if $\mathcal S$ is infinite, one stratum has positive
dimension.
By the definition of a Nash submanifold, it has local Nash coordinate charts
\cite[Def.~2.9.9]{BCR98}.
Restricting such a parametrization of that stratum to a
coordinate line gives a nonconstant Nash arc $\Gamma:J\to\mathcal S$,
a real-analytic semialgebraic map from an open interval $J$.
By \cite[Prop.~8.1.8]{BCR98},
each coordinate of $\Gamma$ is algebraic over $\R(t)$, hence over
$\mathbb{C}(t)$. Define the projection and its image's complex Zariski closure by
\[
 \pi(\delta,u,\mu)=\delta,\qquad
 C=\overline{\pi(\Gamma(J))}^{\,\mathrm{Zar}}
 \subset\Bigl\{\delta\in\mathbb{C}^5:\sum_i m_i\delta_i=0\Bigr\}.
\]
The Zariski closure is the smallest complex algebraic set containing
the projected arc. The projection is nonconstant: on $\mathcal S$,
positivity uniquely determines $u$ from $\delta$, and any row of
\eqref{js:angular} then determines $\mu$. Thus at least one projected
coordinate is nonconstant. At an interior parameter value, the real-analytic
projected coordinate germs $\delta_i(t)$ extend holomorphically. Define
\[
 E=\mathbb{C}\bigl(\delta_1(t),\ldots,\delta_5(t)\bigr)
\]
as a subfield of the meromorphic germ field at that parameter value.
Since each generator is algebraic over $\mathbb{C}(t)$, adjoining $t$ gives
the finite extension $E(t)/\mathbb{C}(t)$, so $\operatorname{trdeg}_{\mathbb{C}}E\le1$;
nonconstancy gives equality. The ideal of polynomials vanishing on
the projected arc is the kernel of the evaluation map into the field of
meromorphic germs at that parameter value: the identity theorem
\cite[Chap.~I, Sec.~3, proof of Thm.~3.3]{Gri89}
extends vanishing on a real subinterval to vanishing as a germ. The image
of this map is the ring generated by the coordinate germs, so its fraction
field is $E$. To verify irreducibility directly, suppose
$C=C_1\cup C_2$ with both $C_k$ proper. Choose polynomials $f_k$ vanishing
on $C_k$ but not identically on $C$. Then $f_1f_2$ vanishes on the arc.
The identity theorem gives $(f_1f_2)\circ\pi\circ\Gamma=0$ as a germ.
Since germs form a domain, one factor vanishes on a subinterval; the
real-analytic identity theorem on connected $J$ makes it vanish on the whole arc,
contradicting that $C_k$ is proper. Thus $C$ is irreducible and
its function field is $E$, of transcendence degree one. Here an algebraic
curve means an irreducible algebraic set whose function field has
transcendence degree one.
In this paper, \emph{Gale curve} denotes this projected curve in the
normalized Gale coordinates.

\paragraph{A smooth projective model.}
Let $K=\mathbb{C}(C)$ be the function field of $C$. Choose a nonconstant $f\in K$.
Since $f$ is nonconstant, it is transcendental over $\mathbb{C}$ and hence a
transcendence basis of $K/\mathbb{C}$. The extension $K/\mathbb{C}(f)$ is algebraic and
finitely generated, so it is finite; it is separable because the
characteristic is zero. The primitive element theorem gives $g\in K$ with
$K=\mathbb{C}(f,g)$
\cite[Chap.~V, Sec.~4, Thm.~4.6]{Lan02}.
Clearing denominators in the minimal polynomial of $g$ over $\mathbb{C}(f)$,
then dividing out the greatest common divisor of its coefficients in $\mathbb{C}[T_1]$,
gives a polynomial relation $\mathcal P(T_1,T_2)=0$. By Gauss's lemma and its
divisibility corollary
\cite[Chap.~IV, Sec.~2, Thm.~2.1 and Cor.~2.2]{Lan02}, $\mathcal P$ is irreducible and the kernel of
$\mathbb{C}[T_1,T_2]\to K$, $T_1\mapsto f$, $T_2\mapsto g$, is $(\mathcal P)$:
every polynomial relation is divisible by $\mathcal P$ in $\mathbb{C}(T_1)[T_2]$, hence
also in $\mathbb{C}[T_1,T_2]$. Thus the plane curve $C'=\{\mathcal P=0\}$ has function
field $K$ and is birational to $C$.
Let $\overline{C'}$ be its projective closure, and let $X$ be its
normalization. The normalization theorem gives a compact connected Riemann surface
$X$ and a holomorphic map $\nu:X\to\overline{C'}$
\cite[Chap.~I, Thm.~2.1 and Def.~2.3; Chap.~II, Sec.~3,
Def.~3.1]{Gri89}. The map is injective over the smooth
locus, and the preimage $\nu^{-1}(\operatorname{Sing}(\overline{C'}))$
of the singular locus is finite
\cite[Chap.~II, Sec.~3, Def.~3.1]{Gri89}.
The curve $C$ may be singular; local analysis takes place on the smooth surface $X$.
Here holomorphic means complex analytic, and meromorphic functions may have
isolated poles \cite[Chap.~I, Sec.~3, Def.~3.1]{Gri89}.
Since $K=\mathbb C(C')$, each element of $K$ pulls back through $\nu$
in a projective coordinate chart to a quotient of holomorphic functions
whose denominator is not identically zero. Factoring their local zero
orders shows that the quotient has at most a pole, so it is meromorphic
on $X$. Meromorphic functions on $X$ form a field
\cite[Chap.~I, Sec.~3, Rem.~3.2]{Gri89}.
Every holomorphic function on a compact connected Riemann surface is constant
\cite[Chap.~I, Sec.~3, Thm.~3.3]{Gri89}.

\paragraph{Descent of the cube roots.}
Write
\[
 \omega=e^{2\pi i/3},\qquad
 \Omega_3=\{z\in\mathbb C:z^3=1\}=\{1,\omega,\omega^2\}.
\]
At a fixed interior point of the real arc, let $u_{ij}$ also denote the
germs of the positive cube-root branches. Write $L=K(u_{ij}:i<j)$ for the
extension generated by these germs. Each radicand $1-\delta_i\delta_j$ is nonzero in $K$, since it is
positive on the arc. The field $K$ contains all cube roots of unity and
has characteristic zero, so $L$ is the common splitting field of the
separable polynomials $T^3-(1-\delta_i\delta_j)$. Hence $L/K$ is finite
Galois \cite[Chap.~V, Thm.~3.3 and Thm.~4.4;
Chap.~VI, Sec.~1]{Lan02}. For
$G=\operatorname{Gal}(L/K)$, the fixed field is
$L^G=\{a\in L:\sigma(a)=a\text{ for every }\sigma\in G\}=K$
\cite[Chap.~VI, Thm.~1.2]{Lan02}. Here descent means that
the germs $u_{ij}$, initially in $L$, belong to $K$.

\begin{lemma}[Descent and two constant identities]
\label{js:constants}
Each $u_{ij}$ belongs to $K$. Moreover, on $X$,
\begin{equation}\label{js:constant-sums}
 \sum_{i<j}m_i m_j u_{ij}
 =\sum_{i<j}m_i m_j u_{ij}^{-2}=U_0,
\end{equation}
where $U_0>W$. The multiplier $\mu$ is meromorphic on $X$.
\end{lemma}
\begin{proof}
Along the chosen positive arc, constrained criticality gives
$d\Phi/dt=0$. Hence $\sum_{i<j}m_i m_ju_{ij}=U_0>0$,
also as an identity of algebraic germs in $L$.
By the setup above, $L/K$ is finite Galois with fixed field $K$.
Let $\sigma\in G=\operatorname{Gal}(L/K)$; it fixes $K$ pointwise and sends
$u_{ij}$ to $\chi_{ij}(\sigma)u_{ij}$, where $\chi_{ij}(\sigma)\in\Omega_3$.
Since $\sigma$ fixes $U_0\in\mathbb{C}\subset K$, applying it to the constant
identity and taking real parts on the original positive arc gives
\[
 \sum_{i<j}m_i m_ju_{ij}(1-\operatorname{Re}\chi_{ij}(\sigma))=0.
\]
Positivity makes every summand nonnegative, so every
$\chi_{ij}(\sigma)=1$. Thus every generator $u_{ij}$ is fixed by every
element of $G$, so $u_{ij}\in L^G=K$ and $L=K$.

Multiply row $i$ of \eqref{js:angular} by $m_i \delta_i$, sum over $i$,
and divide by two. The weighted linear equation and
$\delta_i \delta_j=1-u_{ij}^3$ give
\[
 0=\sum_{i<j}m_i m_j\delta_i \delta_ju_{ij}^{-2}
 =\sum_{i<j}m_i m_j(u_{ij}^{-2}-u_{ij}),
\]
proving the second constant identity.

On the positive arc, the two constant identities give
\begin{align*}
 3(U_0-W)
 &=\sum_{i<j}m_i m_j(2u_{ij}+u_{ij}^{-2}-3)\\
 &=\sum_{i<j}m_i m_j
 \frac{(u_{ij}-1)^2(2u_{ij}+1)}{u_{ij}^2}\ge0.
\end{align*}
Equality would force every $u_{ij}=1$, hence $\delta_i\delta_j=0$
for every pair. At most one $\delta_i$ could then be nonzero, and
$\sum_i m_i\delta_i=0$ would give $\delta=0$, contrary to the choice
of the arc in $\mathcal S$. Thus $U_0>W$.

Finally, \eqref{js:angular} expresses $\mu$ as a rational function of
the $\delta_i$ and the descended $u_{ij}$, so it is meromorphic on $X$.
\end{proof}

Equations \eqref{js:angular} and \eqref{js:constant-sums} hold as meromorphic
identities on $X$, since they hold as identities of germs in $K$.

For each local argument, fix $P\in X$ and a local parameter $s$,
a complex analytic coordinate with $s(P)=0$.
We use $\operatorname{ord}_P(f)$ for the order (or valuation) of a nonzero
meromorphic function or germ $f$ at $P$: the exponent of the first nonzero
power in its local Laurent expansion, positive for a zero and negative
for a pole
\cite[Chap.~I, Sec.~3, Def.~3.4]{Gri89}.
For the zero function or germ, we set $\operatorname{ord}_P(0)=+\infty$.
Writing each nonzero meromorphic germ as a power of $s$ times a
holomorphic germ nonzero at $P$, and comparing leading terms, gives
\[
 \begin{aligned}
 \operatorname{ord}_P(fg)&=\operatorname{ord}_P(f)+\operatorname{ord}_P(g),\\
 \operatorname{ord}_P(f+g)&\ge
 \min\{\operatorname{ord}_P(f),\operatorname{ord}_P(g)\}.
 \end{aligned}
\]
The second inequality is an equality when the two orders differ.
Thus a finite sum with a unique term of lowest order has that same
order. A meromorphic germ is holomorphic at $P$ precisely when its
order is nonnegative, and vanishes there precisely when its order is positive.
Each radical is a nonzero meromorphic function because it is positive
on the original arc. Gale coordinates that vanish identically contribute
zero to the row equations. The local analysis covers every point of $X$,
including nonreal points.
A Gale pole is a point where at least one Gale coordinate has a pole;
coordinates attaining the largest pole order are called dominant.

At a Gale pole $P$, let $\alpha>0$ be the maximal pole order, and let
$d_i\ne0$ be the coefficient of $s^{-\alpha}$ in each dominant coordinate.
For two dominant labels $i<j$, the leading term of
$1-\delta_i\delta_j$ is $-d_i d_j s^{-2\alpha}$. Thus
\[
 3\operatorname{ord}_P(u_{ij})
 =\operatorname{ord}_P(1-\delta_i\delta_j)=-2\alpha.
\]
Since $u_{ij}$ is meromorphic, its order is an integer, so
$3\mid\alpha$ whenever such a pair exists. We use the local expansions
\begin{equation}\label{js:dominant-expansions}
 \begin{gathered}
 \delta_i=d_i s^{-\alpha}+\cdots,\qquad d_i\in\mathbb{C}\setminus\{0\}
 \quad(i\text{ dominant}),\\
 u_{ij}=h_{ij}s^{-2\alpha/3}+\cdots,\qquad
 h_{ij}^3=-d_i d_j\ne0
 \quad(i<j\text{ both dominant}).
 \end{gathered}
\end{equation}
The omitted terms have higher powers of $s$. The coefficients depend
on $P$ and the local parameter. Throughout the local analysis,
$h_{ij}$ denotes the first nonzero Laurent coefficient of $u_{ij}$.

For a holomorphic function $f$
with $f(P)=1$, shrink the neighborhood so that $|f-1|<1$.
Substituting $z=f-1$ in the convergent series
\[
 \log(1+z)=\sum_{k=1}^{\infty}\frac{(-1)^{k+1}z^k}{k}
 =z+O(z^2),\qquad |z|<1,
\]
defines a local holomorphic logarithm with value zero at $P$.
For a holomorphic function $u$ with $u(P)\ne0$,
\[
 u(s)-u(P)=\frac{u(s)^3-u(P)^3}
 {u(s)^2+u(s)u(P)+u(P)^2}.
\]
The denominator has value $3u(P)^2\ne0$ at $P$, so
$u-u(P)$ and $u^3-u(P)^3$ have the same order.

\begin{lemma}[Elementary pole restrictions]
\label{js:pole-count}
At any pole of the Gale vector, at least three coordinates attain
its maximal pole order $\alpha>0$, and not all five do. Every Gale pole
order is divisible by three. At every point of $X$ there is a pair
of coordinates that are not identically zero and whose orders sum to zero. Consequently,
if four coordinates have maximal pole order $\alpha$, the fifth has zero
order exactly $\alpha$.
\end{lemma}
\begin{proof}
Fix a Gale pole $P\in X$. The weighted linear relation forbids a unique
maximal pole, so the preceding order calculation gives $3\mid\alpha$.
Pairing a maximal pole with any other pole proves the remaining
divisibility assertion.
If exactly two coordinates, say $\delta_1,\delta_2$, had maximal pole
order $\alpha$, use \eqref{js:dominant-expansions} for this pair.
For every other pair, at least one factor has order greater than $-\alpha$,
so $\operatorname{ord}_P(\delta_i\delta_j)>-2\alpha$ and consequently
$\operatorname{ord}_P(u_{ij})>-2\alpha/3$.
Thus only $m_1m_2u_{12}$ contributes to the coefficient of $s^{-2\alpha/3}$
in the first sum of \eqref{js:constant-sums}. This coefficient is
$m_1m_2h_{12}\ne0$, whereas the constant $U_0$ has no terms with negative
powers, a contradiction.
If all five were maximal, every inverse square in the second sum would
vanish at the point, contrary to $U_0>0$.

Fix any $P\in X$ and suppose no two coordinates that are not identically zero
have orders summing to zero. A pair whose orders sum to a negative number makes $u_{ij}^{-2}$ tend to
zero. A pair whose orders sum to a positive number makes it tend to a cube root
of unity. If a coordinate is identically zero, $u_{ij}=1$ identically
by its positive branch on the arc. Taking real parts in the second
constant identity would give $U_0\le W$. Thus a pair whose orders sum to zero
exists. When four coordinates have maximal pole order, such a pair consists
of the fifth coordinate and a dominant coordinate, proving the last assertion.
\end{proof}

\subsection{The multiplier at finite Gale points}
\label{js:finite-points}

At a point $P\in X$ where all $\delta_i$ are finite, write
$a_i=\delta_i(P)$ and let $G_P$ be the graph on labels $1,\ldots,5$
whose edges are the pairs $ij$ for which
$u_{ij}(P)=0$, equivalently $\delta_i(P)\delta_j(P)=1$. An isolated
vertex is a label incident to no edge. We will show that such a vertex
always exists. In the diagrams below, only edges of $G_P$ are drawn;
line crossings are not vertices.

We write $H_n$ for the complete graph on $n$ vertices, in which every
pair of distinct vertices is joined. We write $H_{p,q}$ for the complete
bipartite graph with two disjoint vertex classes $I_p,I_q$, where
$p=|I_p|$ and $q=|I_q|$ are positive integers counting their vertices.
Every pair from different classes is joined, and no pair within a class
is joined. In the bipartite cases on five vertices below, $p+q=5$; we order
the classes so that $p\le q$.

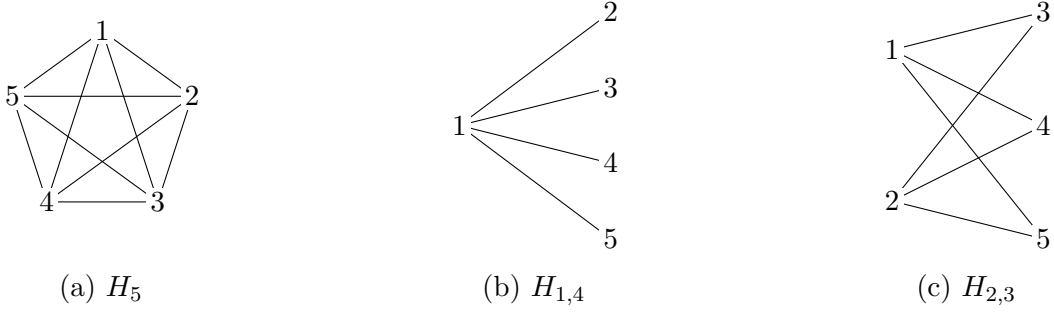
\begin{figure}[htbp]
\centering
\begin{minipage}[b]{0.31\textwidth}
\centering
\begin{tikzpicture}[finite graph]
 \path[use as bounding box] (-1.6,-1.8) rectangle (1.6,1.8);
 \foreach \i in {1,...,5}
  \coordinate (v\i) at ({90-72*(\i-1)}:1.25);
 \foreach \i/\j in {1/2,1/3,1/4,1/5,2/3,2/4,2/5,3/4,3/5,4/5}
  \draw (v\i)--(v\j);
 \foreach \i in {1,...,5} \node at (v\i) {$\i$};
\end{tikzpicture}
\par\small (a) $H_5$
\end{minipage}\hfill
\begin{minipage}[b]{0.31\textwidth}
\centering
\begin{tikzpicture}[finite graph]
 \path[use as bounding box] (-0.4,-1.8) rectangle (2.4,1.8);
 \coordinate (v1) at (0,0);
 \foreach \i/\y in {2/1.5,3/0.5,4/-0.5,5/-1.5}
  \coordinate (v\i) at (2,\y);
 \foreach \j in {2,...,5} \draw (v1)--(v\j);
 \foreach \i in {1,...,5} \node at (v\i) {$\i$};
\end{tikzpicture}
\par\small (b) $H_{1,4}$
\end{minipage}\hfill
\begin{minipage}[b]{0.31\textwidth}
\centering
\begin{tikzpicture}[finite graph]
 \path[use as bounding box] (-0.4,-1.8) rectangle (2.4,1.8);
 \coordinate (v1) at (0,1);
 \coordinate (v2) at (0,-1);
 \foreach \i/\y in {3/1.5,4/0,5/-1.5}
  \coordinate (v\i) at (2,\y);
 \foreach \i in {1,2} \foreach \j in {3,4,5}
  \draw (v\i)--(v\j);
 \foreach \i in {1,...,5} \node at (v\i) {$\i$};
\end{tikzpicture}
\par\small (c) $H_{2,3}$
\end{minipage}
\caption{Connected candidates for $G_P$, up to relabeling.
In (b), $I_p=\{1\}$ and $I_q=\{2,3,4,5\}$;
in (c), $I_p=\{1,2\}$ and $I_q=\{3,4,5\}$.}
\label{js:fig-connected}
\end{figure}

\begin{lemma}[The graph $G_P$ is not connected]
\label{js:connected}
At every finite Gale point $P$, the graph $G_P$ is not connected.
\end{lemma}
\begin{proof}
Suppose $G_P$ were connected. Along an edge a nonzero coordinate value
$a$ is replaced by $a^{-1}$. A cycle of odd length forces $a=a^{-1}$,
so all five values are $1$ or all are $-1$, giving $H_5$.
Otherwise, the values split into $a$ and $a^{-1}$ with $a\ne\pm1$;
every pair from different parts is joined and none within a part.
The only possibilities on five vertices are then $H_{1,4}$ and $H_{2,3}$,
as shown in Figure~\ref{js:fig-connected}.
The complete graph contradicts the weighted linear relation.
In either bipartite case, choose $a$ so that $\delta_i(P)=a$ for
$i\in I_p$ and $\delta_i(P)=a^{-1}$ for $i\in I_q$. Define the total
masses and the pair weights within each part by
\[
 M_p=\sum_{i\in I_p}m_i,\qquad M_q=\sum_{i\in I_q}m_i,\qquad
 W_p=\sum_{\substack{i<j\\i,j\in I_p}}m_i m_j,\qquad
 W_q=\sum_{\substack{i<j\\i,j\in I_q}}m_i m_j.
\]
Thus the weighted linear relation is $M_pa+M_qa^{-1}=0$, giving
$a^2=-M_q/M_p$. All cross radicals vanish, and the internal radicals
have moduli $(1+M_q/M_p)^{1/3}$ and $(1+M_p/M_q)^{1/3}$, respectively.
Hence
\[
 U_0\le W_p(1+M_q/M_p)^{1/3}+W_q(1+M_p/M_q)^{1/3}.
\]
For $x>0$, $(1+x)^{1/3}-1<x/3$, while $W_p<M_p^2/2$ and
$W_q<M_q^2/2$. The right side is therefore less than
$W_p+W_q+M_pM_q=W$, contrary to \eqref{js:constant-sums}.
\end{proof}

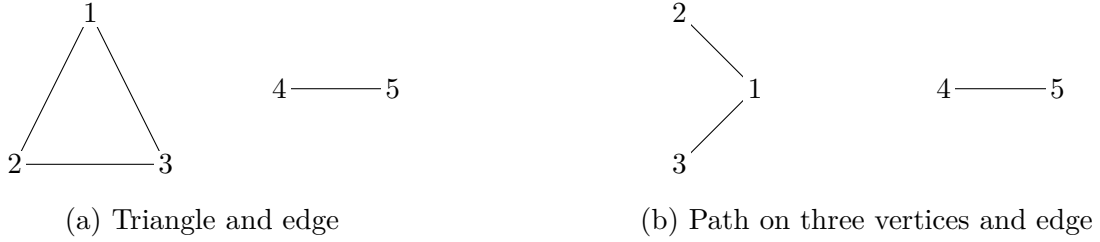
\begin{figure}[htbp]
\centering
\begin{minipage}[b]{0.47\textwidth}
\centering
\begin{tikzpicture}[finite graph]
 \path[use as bounding box] (-0.4,-1.4) rectangle (5.4,1.4);
 \foreach \i/\x/\y in {1/1/1,2/0/-1,3/2/-1,4/3.5/0,5/5/0}
  \coordinate (v\i) at (\x,\y);
 \foreach \i/\j in {1/2,1/3,2/3,4/5} \draw (v\i)--(v\j);
 \foreach \i in {1,...,5} \node at (v\i) {$\i$};
\end{tikzpicture}
\par\small (a) Triangle and edge
\end{minipage}\hfill
\begin{minipage}[b]{0.47\textwidth}
\centering
\begin{tikzpicture}[finite graph]
 \path[use as bounding box] (-0.4,-1.4) rectangle (5.4,1.4);
 \foreach \i/\x/\y in {1/1/0,2/0/1,3/0/-1,4/3.5/0,5/5/0}
  \coordinate (v\i) at (\x,\y);
 \foreach \i/\j in {1/2,1/3,4/5} \draw (v\i)--(v\j);
 \foreach \i in {1,...,5} \node at (v\i) {$\i$};
\end{tikzpicture}
\par\small (b) Path on three vertices and edge
\end{minipage}
\caption{Disconnected candidates without isolated vertices, excluded
by Lemmas~\ref{js:clique-edge} and~\ref{js:path-edge}, respectively.}
\label{js:fig-disconnected}
\end{figure}

If $G_P$ has no isolated vertex, each connected component has at least
two vertices, so its component sizes are five or three plus two.
In the latter case the component on three vertices is a triangle or a path.
The next two lemmas exclude these remaining cases
(Figure~\ref{js:fig-disconnected}).

\begin{lemma}[The graph $G_P$ is not a triangle plus an edge]
\label{js:clique-edge}
For five positive masses and any finite Gale point $P$, the graph $G_P$
cannot consist of a triangle and a separate edge.
\end{lemma}
\begin{proof}
Let $P\in X$ be such a finite Gale point.
After relabeling and a sign change the first three limiting coordinates
are $1$, and the last two are reciprocal. In this proof, let
$I_p=\{4,5\}$ and $I_q=\{1,2,3\}$ be the vertex sets of the two
connected components, so $p=2\le q=3$. Put
\[
 \begin{aligned}
 M_p&=m_4+m_5, & M_q&=m_1+m_2+m_3,\\
 W_p&=m_4m_5, & W_q&=m_1m_2+m_1m_3+m_2m_3.
 \end{aligned}
\]
The vanishing radicals are exactly the three edges among labels
$1,2,3$ and the edge $45$.
Let $w>0$ be their largest zero order, and retain only edges attaining
that order. For such an edge write
$u_{ij}=h_{ij}s^w+o(s^w)$, with $h_{ij}\ne0$.
Let $\gamma$ denote the coefficient of $s^{-2w}$ in the common
multiplier. Every row of the critical-point equations has coefficient
\[
 \sum_{j:ij\text{ maximal}}m_j a_jh_{ij}^{-2}=\gamma.
\]
Edge $45$ must attain the maximal order. Otherwise its rows give
$\gamma=0$. A nonempty maximal graph on the first three vertices
with all row sums zero can have no vertex of degree one, so would have
to be a triangle. A symmetric matrix of size three with zero diagonal
and nonzero off-diagonal entries $q_{12},q_{13},q_{23}$ has determinant
$2q_{12}q_{13}q_{23}\ne0$,
and cannot annihilate the positive vector $(m_1,m_2,m_3)$. This is
impossible. The rows of edge $45$ consequently give $m_4a_4=m_5a_5$.
Together with $a_4a_5=1$ and the weighted linear relation, this yields
\[
 m_4a_4=m_5a_5=-M_q/2,\qquad
 m_4m_5=M_q^2/4,\qquad \gamma\ne0,
 \qquad h_{45}^2=-M_q/(2\gamma).
\]
Each of the first three vertices meets an edge of maximal order,
so their maximal graph is either a path or a triangle.

For each maximal edge among labels $1,2,3$, put
$\rho_{ij}=m_i m_jh_{ij}^{-2}$.
In the path case, label the maximal edges $12,13$, with middle vertex $1$.
The end rows give $\rho_{12}=\gamma m_2$ and
$\rho_{13}=\gamma m_3$; the middle row gives
$\rho_{12}+\rho_{13}=\gamma m_1$. Hence $m_1=m_2+m_3=M_q/2$ and
\[
 W_q\ge M_q^2/4.
\]

Suppose next that all three triangle edges have order $w$. Solving
the three row equations yields
\[
 \rho_{ij}=\frac\gamma2(m_i+m_j-m_k),\qquad
 h_{ij}^2=\frac{2m_i m_j}{\gamma(M_q-2m_k)}.
\]
In particular no $m_i$ equals $M_q/2$.
The coefficient of $s^w$ in the constant potential identity gives
\begin{equation}\label{js:finite-clique-potential}
 \sum_{1\le i<j\le3}m_i m_jh_{ij}+\frac{M_q^2}{4}h_{45}=0.
\end{equation}
To justify that no nonvanishing radical contributes to this order,
put $\ell_i=\log(\delta_i/a_i)$ with $\ell_i(P)=0$, using the
normalized local logarithms described above. For $1\le i<j\le3$,
$\ell_i+\ell_j=\log(1-u_{ij}^3)$ has order $3w$, so each of
$\ell_1,\ell_2,\ell_3$ has order at least $3w$.
Likewise $\ell_4+\ell_5$ has order $3w$. Write
$\ell=(\ell_4+\ell_5)/2$ and $v=(\ell_4-\ell_5)/2$.
The weighted linear identity gives $e^\ell\cosh v=1+O(s^{3w})$.
Since $\ell=O(s^{3w})$, this implies
\[
 \cosh v-1=\frac12v^2+O(v^4)=O(s^{3w}),
\]
and hence $\operatorname{ord}_P(v)\ge3w/2$, also when $v\equiv0$.
Every coordinate therefore changes by $O(s^{3w/2})$. The cross
radicands, and hence their nonvanishing cube roots by the local identity
above, also change by $O(s^{3w/2})$.
Since $3w/2>w$, the constant potential identity now proves
\eqref{js:finite-clique-potential}.

If all three masses are less than $M_q/2$, then
$\sum_{i=1}^3m_i^2<M_q^2/2$, so $W_q>M_q^2/4$.
Otherwise, choose $\{i,j,k\}=I_q$ so that $m_k$ is the largest mass;
then $m_k>m_i+m_j$. Choose a square root of $\gamma$. The quantities
$\sqrt\gamma h_{ik}$ and $\sqrt\gamma h_{jk}$ are nonzero real
numbers, whereas $\sqrt\gamma h_{ij}$ and
$\sqrt\gamma h_{45}$ are nonzero purely imaginary numbers.
After multiplying \eqref{js:finite-clique-potential} by
$\sqrt\gamma$, its real part forces
\[
 \frac{m_i^3}{m_k+m_i-m_j}=\frac{m_j^3}{m_k+m_j-m_i}.
\]
If $m_i\ne m_j$, factoring gives
\[
 m_k=\frac{(m_i+m_j)(m_i^2-m_i m_j+m_j^2)}
 {m_i^2+m_i m_j+m_j^2}<m_i+m_j,
\]
a contradiction. Thus $m_i=m_j$. The imaginary part, followed by
squaring positive magnitudes, gives
\[
 64m_i^6=(m_k-2m_i)(m_k+2m_i)^5.
\]
For $t=m_k/(2m_i)>1$ this says $(t-1)(t+1)^5=1$, so $t<2$.
Consequently $m_k<4m_i$ and
\[
 W_q=2m_i m_k+m_i^2>\frac{(m_k+2m_i)^2}{4}=M_q^2/4.
\]
Thus $W_q\ge M_q^2/4$ in all cases, while $W_p=M_q^2/4$.

At the finite point only cross radicals contribute to the potential.
Put $\tau=\sqrt{m_5/m_4}>0$, so $a_4=-\tau$, $a_5=-\tau^{-1}$.
Using $(1+x)^{1/3}<1+x/3$ for $x>0$, their total contribution has
real part at most
\begin{align*}
 U_0&\le M_q\left[m_4(1+\tau)^{1/3}+m_5(1+\tau^{-1})^{1/3}\right]\\
 &<M_pM_q+\frac{M_q}{3}(m_4\tau+m_5/\tau)\\
 &=M_pM_q+M_q^2/3\\
 &<M_pM_q+W_p+W_q=W.
\end{align*}
This contradicts $U_0>W$.
\end{proof}

\begin{lemma}[The graph $G_P$ is not a path plus an edge]
\label{js:path-edge}
For five positive masses and any finite Gale point $P$, the graph $G_P$
cannot consist of a path on three vertices and a separate edge.
\end{lemma}
\begin{proof}
Label the path edges $12,13$ and the separate edge $45$.
Let $w$ be the largest zero order among the vanishing radicals and
$\gamma$ the coefficient of $s^{-2w}$ in the multiplier. The separate edge must attain
that order: otherwise $\gamma=0$, and a maximal edge in the path
would have an endpoint of degree one whose row coefficient could not cancel.
Thus $\gamma\ne0$, and both path edges must attain that order as well.
Let $h_{ij}$ denote their nonzero leading coefficients.
The path's limiting coordinates are $a_1=a$, $a_2=a_3=a^{-1}$.
Its leaf rows give
$m_1a/h_{12}^2=m_1a/h_{13}^2=\gamma$; its middle row therefore gives
$a^2=(m_2+m_3)/m_1$. Similarly the pair rows give $a_4^2=m_5/m_4$.
An overall sign change lets us take $a>0$. Since $a_4a_5=1$,
the weighted linear identity forces $a_4<0$ and becomes
\[
 2\sqrt{m_1(m_2+m_3)}-2\sqrt{m_4m_5}=0.
\]
Hence $m_1(m_2+m_3)=m_4m_5$. Define the positive auxiliary quantity
\[
 d=\sqrt{m_1(m_2+m_3)}=\sqrt{m_4m_5}>0.
\]
The limiting coordinates are then
\[
 \sqrt{\frac{m_2+m_3}{m_1}},\quad
 \sqrt{\frac{m_1}{m_2+m_3}},\quad
 \sqrt{\frac{m_1}{m_2+m_3}},\quad
 -\sqrt{\frac{m_5}{m_4}},\quad-\sqrt{\frac{m_4}{m_5}}.
\]
Between the two components all products are negative. The weighted sums
of coordinate magnitudes on the two components both equal $2d$:
\[
 \sum_{i=1}^3 m_i|a_i|=m_4|a_4|+m_5|a_5|=2d.
\]
By $(1+x)^{1/3}<1+x/3$, the total contribution from pairs in different
components has real part less than
\[
 W_{\mathrm{cross}}+\frac{4d^2}{3},
 \qquad W_{\mathrm{cross}}=(m_1+m_2+m_3)(m_4+m_5).
\]
Among pairs within the components, only $u_{23}$ can contribute.
Its weight is $m_2m_3\le (m_2+m_3)^2/4$ and its modulus is
$|1-m_1/(m_2+m_3)|^{1/3}$. For every $x>0$,
\[
 |1-x|^{1/3}-1\le x/3.
\]
For $0<x\le2$ the left side is nonpositive. For $x>2$, applying
$(1+t)^{1/3}-1<t/3$ with $t=x-2>0$ gives
\[
 (x-1)^{1/3}-1<(x-2)/3<x/3.
\]
Hence its excess over its pair weight is at most
\[
 m_2m_3\bigl(|1-m_1/(m_2+m_3)|^{1/3}-1\bigr)
 \le \frac{m_1m_2m_3}{3(m_2+m_3)}
 \le\frac{d^2}{12}.
\]
The zero edges have total pair weight $m_1(m_2+m_3)+m_4m_5=2d^2$.
Consequently
\[
 U_0<W_{\mathrm{cross}}+m_2m_3+\frac{17d^2}{12}
 <W_{\mathrm{cross}}+m_2m_3+2d^2=W,
\]
a contradiction.
\end{proof}

\begin{corollary}[The multiplier at finite Gale points is holomorphic]
\label{js:finite-regular}
At every finite Gale point, $\mu$ is holomorphic.
\end{corollary}
\begin{proof}
Fix a finite Gale point $P$. By the preceding classification and
Lemmas~\ref{js:connected}, \ref{js:clique-edge} and~\ref{js:path-edge},
$G_P$ has an isolated vertex $i$. Its row equation is
\[
 \mu=\sum_{j\ne i}m_j\delta_j u_{ij}^{-2}.
\]
Each $\delta_j$ is holomorphic at $P$ and each $u_{ij}(P)\ne0$,
so this expression is holomorphic there.
\end{proof}

\subsection{The multiplier at Gale poles}
\label{js:poles}

\begin{lemma}[Poles with three dominant coordinates]
\label{js:triple}
At a Gale pole, suppose exactly three coordinates attain the maximal
pole order $\alpha>0$. Then the common multiplier has order exactly
$\alpha/3$. In particular it is holomorphic and vanishes at that point.
The other coordinates may have poles of smaller order or be finite,
and radicals may vanish.
\end{lemma}
\begin{proof}
Fix such a pole $P\in X$.
Label the dominant coordinates $1,2,3$, and use
\eqref{js:dominant-expansions} for these coordinates and their internal radicals.
For each outside coordinate $\delta_j$, $j=4,5$, that is not identically zero,
put $v_j=\operatorname{ord}_P(\delta_j)>-\alpha$.
A vanishing radical incident to the triple must join a dominant
coordinate to an outside coordinate of zero order $\alpha$; for such
an outside coordinate write $\delta_j=d_j s^\alpha+\cdots$.
Every other outside term in a dominant row of the critical-point
equations has order strictly greater than $\alpha/3$: for
$v_j<\alpha$ its order is $(v_j+2\alpha)/3>\alpha/3$, and for
$v_j\ge\alpha$ with nonvanishing radical its order is
$v_j\ge\alpha$. An outside coordinate that is identically zero contributes no term.

First suppose there is no vanishing radical incident to the triple.
Let $\gamma$ be the coefficient of $s^{\alpha/3}$ in $\mu$.
The coefficients at this order in the three critical-point equations give
\[
 \sum_{j\in\{1,2,3\}\setminus\{i\}}
 \frac{m_j d_j}{h_{ij}^2}=\gamma,
 \qquad \sum_{i=1}^3m_i d_i=0.
\]
Multiplying row $i$ by $m_i d_i$ and using $h_{ij}^3=-d_i d_j$ gives
\[
 m_i m_j h_{ij}+m_i m_k h_{ik}=-\gamma m_i d_i,
 \qquad \{i,j,k\}=\{1,2,3\}.
\]
Together with the weighted linear relation, these equations yield
\[
 m_1m_2h_{12}=\gamma m_3d_3,\qquad
 m_1m_3h_{13}=\gamma m_2d_2,\qquad
 m_2m_3h_{23}=\gamma m_1d_1.
\]
Their left sides are nonzero, so $\gamma\ne0$, proving the exact order
$\alpha/3$.

Otherwise let $w>0$ be the largest order among the incident radical
zeros. Every dominant row of the critical-point equations has order at least
$\min\{\alpha/3,\alpha-2w\}$. Choose a matched outside coordinate $\delta_j$
incident to a zero of order $w$. In its row of the critical-point equations the corresponding
term has order $-\alpha-2w<\min\{\alpha/3,\alpha-2w\}$.
Its terms from other outside coordinates
have order greater than $-\alpha$: their products with $\delta_j$ tend to zero,
so their radicals do not vanish. Dominant terms whose radicals do not
vanish have order $-\alpha$. Consequently at least two dominant terms of order
$-\alpha-2w$ must cancel. Their product limits with $\delta_j$ equal one, so
after relabeling
\[
 d_1=d_2,\qquad
 d_3=-\frac{m_1+m_2}{m_3}d_1\ne d_1.
\]

No radical incident to label $3$ can vanish. If $u_{3j}$ had zero
order $w'>0$, then in the outside row $j$ it would be the only
dominant term whose radical vanishes, since its limiting product condition
$d_3d_j=1$ excludes vanishing radicals joining label $j$ to label 1 or 2.
All other terms in that row have order at least $-\alpha$.
Thus $\operatorname{ord}_P(\mu)=-\alpha-2w'$. Every radical zero incident
to label $3$ must have this same order $w'$, by applying the same
argument to its outside row. But the third dominant row has order
at least $\min\{\alpha/3,\alpha-2w'\}>-\alpha-2w'$, a contradiction.

The third dominant row therefore has no incident radical zero.
Its internal coefficient at order $\alpha/3$ is
\[
 \gamma_3=-\frac{m_1h_{13}+m_2h_{23}}{d_3}.
\]
Since $d_1=d_2$, we have $h_{23}/h_{13}\in\Omega_3$.
If $m_1h_{13}+m_2h_{23}=0$, this ratio would equal the negative
real number $-m_1/m_2$, which is not in $\Omega_3$.
Hence $\gamma_3\ne0$, and all other terms in the row have larger
order. Hence $\operatorname{ord}_P(\mu)=\alpha/3$.
\end{proof}

\begin{lemma}[The four internal sums do not all vanish]
\label{js:phase-four}
Let $m_1,\ldots,m_4>0$, let $d_1,\ldots,d_4\ne0$ satisfy
$\sum_{i=1}^4m_i d_i=0$, and let $h_{ij}^3=-d_i d_j$
for $1\le i<j\le4$.
The four internal sums in the critical-point equations
$\sum_{1\le j\le4,\ j\ne i}m_j d_j/h_{ij}^2$ cannot all be zero.
\end{lemma}
\begin{proof}
Choose $b_i$ with $b_i^3=d_i$ for $1\le i\le4$ and write
$h_{ij}=-\varepsilon_{ij}b_ib_j$, with $\varepsilon_{ij}\in\Omega_3$.
Zero sums in the critical-point equations would give
\[
 Q(m_1b_1,\ldots,m_4b_4)^T=0,\qquad \sum_{i=1}^4m_i b_i^3=0,
\]
where $Q$ is symmetric, has zero diagonal, and has off-diagonal
entries $\varepsilon_{ij}$. A diagonal congruence rescales each row and
its corresponding column by the same nonzero factor. Set
\[
 D=\operatorname{diag}(1,\varepsilon_{12}^{-1},
 \varepsilon_{13}^{-1},\varepsilon_{14}^{-1}),\qquad
 \widetilde Q=DQD.
\]
Then $\widetilde QD^{-1}(m_1b_1,\ldots,m_4b_4)^T=0$.
Every diagonal entry of $D$ has cube one, so the transformed kernel
coordinates have the same cubes as the original ones.
The off-diagonal entries in the first row of $\widetilde Q$ are all one.
Writing the remaining entries as $\widetilde{\varepsilon}_{ij}\in\Omega_3$ gives
\[
 \widetilde Q=
 \begin{pmatrix}
 0&1&1&1\\
 1&0&\widetilde{\varepsilon}_{23}&\widetilde{\varepsilon}_{24}\\
 1&\widetilde{\varepsilon}_{23}&0&\widetilde{\varepsilon}_{34}\\
 1&\widetilde{\varepsilon}_{24}&\widetilde{\varepsilon}_{34}&0
 \end{pmatrix}.
\]
Its determinant is
\[
 \begin{aligned}
 \det\widetilde Q
 &=\widetilde{\varepsilon}_{23}^2+\widetilde{\varepsilon}_{24}^2
   +\widetilde{\varepsilon}_{34}^2\\
 &\quad-2(\widetilde{\varepsilon}_{23}\widetilde{\varepsilon}_{24}
   +\widetilde{\varepsilon}_{23}\widetilde{\varepsilon}_{34}
   +\widetilde{\varepsilon}_{24}\widetilde{\varepsilon}_{34}).
 \end{aligned}
\]
If these three entries coincide, the determinant is
$-3\widetilde{\varepsilon}_{23}^2\ne0$.
If exactly two coincide, say
$\widetilde{\varepsilon}_{23}=\widetilde{\varepsilon}_{24}$, it is
$\widetilde{\varepsilon}_{34}(\widetilde{\varepsilon}_{34}
-4\widetilde{\varepsilon}_{23})\ne0$.
If all three are distinct, permute the last three labels to obtain
\[
 (\widetilde{\varepsilon}_{23},\widetilde{\varepsilon}_{24},
 \widetilde{\varepsilon}_{34})=(1,\omega,\omega^2).
\]
The leading principal $3\times3$ minor and a kernel vector satisfy
\[
 \det\begin{pmatrix}0&1&1\\1&0&1\\1&1&0\end{pmatrix}=2,\qquad
 \widetilde Q(1,\omega,\omega^2,1)^T=0.
\]
Hence the matrix has rank three, with kernel spanned by
$(1,\omega,\omega^2,1)^T$. Thus every nonzero kernel vector has all
four cubes equal and nonzero. Undoing the phase change
preserves this property. The weighted linear relation becomes
\[
 0=\sum_{i=1}^4m_i b_i^3=(m_1b_1)^3\sum_{i=1}^4m_i^{-2},
\]
impossible for positive masses.
\end{proof}

\begin{lemma}[Two zero rows at a repeated pair are impossible]
\label{js:two-zero-rows}
Let $m_1,\ldots,m_4>0$. Suppose $d_1=d_2\ne0$,
$d_3,d_4\ne0$, $\sum_{i=1}^4m_i d_i=0$, $h_{ij}^3=-d_i d_j$ for
$1\le i<j\le4$, and $\sum_{1\le i<j\le4}m_i m_jh_{ij}=0$.
Then the internal sums in the critical-point equations at
labels 3 and 4 cannot both be zero.
\end{lemma}
\begin{proof}
Rescale $d_i$ and $h_{ij}$ compatibly so that $d_1=d_2=1$.
Choose $b_3^3=d_3$, $b_4^3=d_4$, set $b_1=b_2=1$, and write
$h_{ij}=-\varepsilon_{ij}b_ib_j$ with $\varepsilon_{ij}\in\Omega_3$.
Put
\[
 F_3=m_1\varepsilon_{13}+m_2\varepsilon_{23},\qquad
 F_4=m_1\varepsilon_{14}+m_2\varepsilon_{24}.
\]
The identity for the vanishing internal potential is
\[
 -m_1m_2\varepsilon_{12}-m_3F_3b_3-m_4F_4b_4
 -m_3m_4\varepsilon_{34}b_3b_4=0.
\]
The two zero rows give $F_3=-m_4\varepsilon_{34}b_4$ and
$F_4=-m_3\varepsilon_{34}b_3$. Substitution yields
$m_1m_2\varepsilon_{12}=m_3m_4\varepsilon_{34}b_3b_4$, hence
\[
 F_3F_4=m_3m_4\varepsilon_{34}^2b_3b_4
 =m_1m_2\varepsilon_{12}\varepsilon_{34}.
\]
In particular, $|F_3F_4|=m_1m_2$.
For two positive weights $m_1,m_2$, the squared modulus of a sum of two
cube roots of unity with these weights is either $(m_1+m_2)^2$ or
$m_1^2-m_1m_2+m_2^2$. Both are at least $m_1m_2$, and equality occurs only
when $m_1=m_2$ and the phases differ.
Since $|F_3|^2|F_4|^2=(m_1m_2)^2$, both lower bounds are attained.
Two distinct elements of $\Omega_3$ sum to the negative of the third, so
\[
 m_1=m_2,\qquad F_3^3=F_4^3=-m_1^3.
\]
The row equations then give $b_4^3=(m_1/m_4)^3>0$ and
$b_3^3=(m_1/m_3)^3>0$, contradicting
$2m_1+m_3b_3^3+m_4b_4^3=0$.
\end{proof}

\begin{lemma}[Poles with four dominant coordinates]
\label{js:four-pole}
At a Gale pole with four coordinates of maximal pole order $\alpha>0$,
the multiplier $\mu$ is holomorphic and vanishes at that point. Its local equations
forbid $\mu$ from being identically zero on $X$.
\end{lemma}
\begin{proof}
Fix such a Gale pole $P\in X$.
Relabel the dominant set as $I=\{1,2,3,4\}$, and use
\eqref{js:dominant-expansions} on $I$. By Lemma~\ref{js:pole-count}
and the weighted linear relation,
\[
 \delta_5=d_5 s^\alpha+\cdots,\qquad d_5\ne0,\qquad
 \sum_{i\in I}m_i d_i=0.
\]
A matched radical
$u_{i5}$ vanishes exactly when $d_i d_5=1$. Not all four can vanish,
since then the weighted linear relation would be nonzero. A dominant row with
nonvanishing matched radical has internal terms of order $\alpha/3$
and outside term of order $\alpha$. Consequently
$\operatorname{ord}_P(\mu)\ge\alpha/3$; thus $\mu$ is holomorphic and vanishes at $P$.

Put $M_I=\sum_{i\in I}m_i$ and define
\[
 \mathcal U_I=\sum_{\substack{i<j\\i,j\in I}}m_i m_j u_{ij},\qquad
 \mathcal I_I=\frac1{M_I}\sum_{\substack{i<j\\i,j\in I}}m_i m_j u_{ij}^{-2}.
\]
For a physical configuration, $\mathcal I_I$ is the moment of inertia
of the bodies in $I$ about their own center of mass.
Multiplying the full rows indexed by $i\in I$ of the critical-point
equations by $m_i\delta_i$ and summing
gives, using the fifth row and the weighted linear relation,
\[
 \begin{aligned}
 -m_5\delta_5\mu&=2(M_I\mathcal I_I-\mathcal U_I)+m_5\delta_5\mu,\\
 \mathcal U_I&=M_I\mathcal I_I+m_5\delta_5\mu.
 \end{aligned}
\]
By \eqref{js:dominant-expansions}, every internal radical has pole order
$2\alpha/3$.
Thus $\mathcal I_I\to0$ and $\delta_5\mu\to0$, so $\mathcal U_I\to0$.
In particular its leading internal coefficient satisfies
\[
 \sum_{\substack{i<j\\i,j\in I}}m_i m_jh_{ij}=0.
\]

There cannot be exactly one zero matched radical: its inverse
square would be the sole pole in the second constant sum. Nor
can there be exactly three. In that case, say $u_{15},u_{25},u_{35}$
vanish, their leading products give $d_1=d_2=d_3=1/d_5$.
The weighted linear relation gives
$d_4d_5=-(M_I-m_4)/m_4$. The remaining radical has a nonzero limit
$h_{45}=u_{45}(P)$ with $h_{45}^3=M_I/m_4$. Because $\mathcal U_I\to0$, the first constant
sum gives $U_0=m_4m_5h_{45}$, whence
\[
 U_0=m_4m_5|h_{45}|=m_5m_4^{2/3}M_I^{1/3}
 <m_5M_I<W.
\]
Here $m_4<M_I$, and
$W=m_5M_I+\sum_{i<j,\ i,j\in I}m_i m_j>m_5M_I$ because all
internal pair weights are positive. This is a contradiction.
Exactly zero or two matched radicals therefore vanish.

For any $i\in I$ with $u_{i5}(P)\ne0$, its row has the expansion
\[
 \mu=s^{\alpha/3}\left(
 \sum_{\substack{j\in I\\j\ne i}}\frac{m_jd_j}{h_{ij}^2}+o(1)
 \right).
\]
The outside term has order $\alpha>\alpha/3$ and does not affect
the displayed coefficient. Suppose now that $\mu\equiv0$, so this
coefficient is zero.

If none vanishes, this applies to all four rows, contrary to
Lemma~\ref{js:phase-four}.
If two vanish, relabel them $u_{15},u_{25}$; their product limits
give $d_1=d_2$. Rows 3 and 4 have nonvanishing matched radicals,
so their internal leading sums both vanish.
Together with the zero leading internal potential, this contradicts
Lemma~\ref{js:two-zero-rows}. Thus $\mu$ cannot be identically zero.
\end{proof}

\subsection{Completion of the proof}
\label{js:compact-contradiction}

\begin{proof}[Proof of Theorem~\ref{js:theorem}]
If $\mathcal S$ were infinite, the construction above would produce
a nonconstant algebraic Gale curve and its connected smooth projective
model $X$. Some Gale coordinate is nonconstant.
If it had no pole, it would be holomorphic on all of $X$, hence constant
by the compactness result recalled above, a contradiction.

At every point with all Gale coordinates finite, the multiplier
$\mu$ is holomorphic by Corollary~\ref{js:finite-regular}. At every
other point, Lemma~\ref{js:pole-count} leaves exactly three or
four maximal pole coordinates. Lemmas~\ref{js:triple} and
\ref{js:four-pole} show that $\mu$ is holomorphic and vanishes there; they also
forbid $\mu$ from being identically zero. Consequently $\mu$ is
globally holomorphic on $X$, so the same compactness result makes it constant.
At a Gale pole its value is zero,
so that constant is zero, contradicting either local pole lemma.
This contradiction proves that $\mathcal S$ is finite. Every normalized
physical configuration gives a point of $\mathcal S$, and each such point
determines all distances $r_{ij}$.
By the argument using inner products in Subsection~\ref{js:setup}, these determine
the similarity class after the chosen normalization, so the physical
similarity classes are finite in number.
\end{proof}

\section*{Acknowledgements}
\begin{sloppypar}
The first named author is partially supported by the National Key Research and Development
Program of China (NKPs) [Grant Number 2024YFA1013201], the Natural Science Foundation of
Chongqing (NSFCQ) [Grant Number \mbox{CSTB2024NSCQ-LZX0040}], and the Special Project of Chongqing
Municipal Science and Technology Bureau [Grant Number 2025CCZ015].

\par\bigskip
\noindent\textbf{Declaration of AI use}\par\nobreak\smallskip

This manuscript was prepared and revised with the assistance of OpenAI's Codex,
including help with proof analysis, mathematical exposition, notation,
reference checking, and typesetting. The authors take full responsibility
for the content.
\end{sloppypar}

\end{document}